\documentclass[12pt]{article}
\usepackage{amsmath,amssymb,amsfonts,amsthm}
\usepackage{eucal} 
\usepackage{pstricks,pst-plot,pst-node}
\newpsobject{showgrid}{psgrid}{subgriddiv=1,griddots=5,gridlabels=0pt}

\newcommand{\lang}{\left\langle}
\newcommand{\rang}{\right\rangle}
\newcommand{\blang}{\big\langle}
\newcommand{\brang}{\big\rangle}
\newcommand{\dlang}{{\mathbf {\Big (}}}
\newcommand{\drang}{{\mathbf {\Big )}}}

\newcommand{\Bv}{\Big |}
\newcommand{\bv}{\big |}
\newcommand{\lv}{\left |}

\newcommand{\zz}{{\mathfrak{z}}}
\newcommand{\vac}{v_\emptyset}
\newcommand{\bH}{\mathsf{H}}

\newcommand{\bJ}{\mathsf{J}}

\newcommand{\cF}{\mathcal{F}}

\newcommand{\notdiv}{\!\!\not|\,}

\newcommand{\MM}{\mathsf{M}}
\newcommand{\tMM}{\MM_D}

\newcommand{\C}{\mathbb{C}}
\newcommand{\Q}{\mathbb{Q}}
\newcommand{\Z}{\mathbb{Z}}
\newcommand{\cO}{\mathcal{O}}

\newcommand{\Pp}{\mathbf{P}^1}

\newcommand{\cI}{\mathcal{I}}
\newcommand{\ld}{\,\text{---}\,}
\newcommand{\TT}{T^\pm}

\DeclareMathOperator{\Hilb}{Hilb}
\DeclareMathOperator{\sgn}{sgn}
\DeclareMathOperator{\Sym}{Sym}

\DeclareMathOperator{\Aut}{Aut}

\newcommand{\M}{\overline{M}^{\bullet}_{h}}

\newtheorem{Theorem}{Theorem}
\newtheorem{Lemma}{Lemma}
\newtheorem{Corollary}[Lemma]{Corollary}

\begin{document}
\title{Quantum cohomology of the Hilbert scheme of 
points in the plane}
\author{A. Okounkov and R. Pandharipande} 
\date{}
 \maketitle

\begin{abstract}
We determine the ring structure of the 
equivariant quantum cohomology of the Hilbert scheme of 
points of $\C^2$. The operator of quantum multiplication 
by the divisor class is a 
nonstationary deformation of the quantum 
Calogero-Sutherland many-body system. 
A relationship between the quantum cohomology of the Hilbert scheme 
and the Gromov-Witten/Donaldson-Thomas correspondence
for local curves is proven.
\end{abstract}

\tableofcontents

\section{Introduction}

\subsection{Overview}

The Hilbert scheme $\Hilb_n$ 
of $n$ points in the plane $\C^2$ 
parametrizes ideals $\cI\subset \C[x,y]$ of colength $n$,
$$
\dim_\C \C[x,y]/\cI = n \,. 
$$
An open dense set of $\Hilb_n$ parameterizes ideals associated to configurations of
$n$ distinct points.
The Hilbert scheme is a nonsingular, irreducible, quasi-projective algebraic variety of 
dimension $2n$ with a rich and much studied geometry,
see
\cite{goe,Nak} for an introduction.

The symmetries of $\C^2$ lift to the Hilbert scheme. 
The algebraic torus 
$$T=(\C^*)^2$$ 
acts on $\C^2$ by scaling coordinates,
$$
(z_1,z_2) \cdot (x,y) = (z_1 x, z_2 y)\,.
$$
The induced $T$-action on $\Hilb_n$ plays a central role
in the subject.

The $T$-equivariant cohomology of $\Hilb_n$  has been recently determined, see
\cite{CG,Lehn,LS,LQW,vass}. As a ring, 
$H^*_T(\Hilb_n, {\mathbb Q})$ is generated by the Chern 
classes of the tautological rank $n$ bundle 
$$
\cO/\cI\rightarrow \Hilb_n,
$$
with fiber
$\C[x,y]/\cI$ over $[\cI]\in \Hilb_n$.  
The operator of classical multiplication by the 
divisor 
$$
D = c_1 (\cO/\cI) 
$$
in 
$H^*_T(\Hilb_n, {\mathbb Q})$
is naturally identified with the Hamiltonian of the 
Calogero-Sutherland integrable quantum many-body system. 
The ring $H^*_T(\Hilb_n, {\mathbb Q})$ is 
a module over
$$H_T^*({\text {pt}}) = {\mathbb Q}[t_1,t_2],$$
where $t_1$ and $t_2$ are the Chern classes of the
respective factors of the standard representation.
The ratio  $-t_2/t_1$
of the equivariant parameters is identified with 
the coupling constant in the Calogero-Sutherland system. 

The goal of our paper is to 
compute the small quantum product on the 
$T$-equivariant cohomology of $\Hilb_n$.  
The matrix elements of the small quantum product 
count, in an appropriate sense, rational 
curves meeting three given 
subvarieties of $\Hilb_n$. 
The (non-negative) degree of a curve class $\beta\in H_2(\Hilb_n,{\mathbb Z})$ is defined by
 $$d= \int_\beta D.$$
Curves of degree $d$ are counted with 
weight $q^d$, where $q$ is the quantum parameter.
The ordinary multiplication in 
$T$-equivariant cohomology is recovered by
setting $q=0$. See \cite{CK,fp} for an introduction to 
quantum cohomology.

Our main result, Theorem \ref{T1}, is an explicit formula for the operator 
$\MM_D$ of small quantum multiplication by $D$.  As a corollary, 
$D$ is proven to generate the small quantum cohomology ring 
over $\Q(q,t_1,t_2)$. The ring structure
is therefore determined. 

The full $T$-equivariant quantum cohomology in genus 0 (with arbitrary numbers
of insertions) is easily calculated from the 3-point
invariants. A procedure is presented in Section \ref{s_big}. 
The higher genus invariants are discussed in Section \ref{higherg}.

\subsection{Quantum differential equation}

Our explicit form for $\MM_D$ implies, 
the {\em  quantum differential 
equation} 
\begin{equation}
  \label{qode}
  q\frac{d}{dq} \, \psi = \MM_D \, \psi\,, \quad \psi(q) \in H^*_T(\Hilb_n,
{\mathbb Q})\,, 
\end{equation}
has regular singularities at 
$q=0,\infty,$ and certain roots of unity. 
The monodromy of this linear ODE is 
remarkable \cite{s_QDE}. In particular, 
we prove the monodromy is invariant under
$$
t_1 \mapsto t_1 - 1
$$
provided 
$$
t_1 \ne \frac{r}{s}\,, \quad 0<r\le s\le n \,, 
$$
and similarly for $t_2$.  

As a corollary, when the sum 
$t_1+t_2$ of equivariant parameters is an integer, 
there is no monodromy around the roots of unity.    
The  full monodromy is then abelian (and, in fact, diagonalizable
for generic $t_1$). 

Equation \eqref{qode} may be viewed as an 
exactly solvable nonstationary generalization of the 
Calogero-Sutherland system. Several results and 
conjectures concerning its solutions, which are 
deformations of Jack polynomials, are presented in 
\cite{s_QDE}. 


We expect to find similar integrability in 
the quantum differential equation for the Hilbert scheme 
of points of any smooth surface.

\subsection{Relation to Gromov-Witten and Donaldson-Thomas theories}

Consider the projective line $\Pp$ with three
distinguished points 
$$0,1,\infty \in \Pp.$$
The $T$-equivariant  Gromov-Witten theory of $\Pp\times \C^2$ 
relative to $\{0,1,\infty\}$. 
has been calculated in \cite{jbrp}. 
Let relative conditions be specified by 
$$\lambda,\mu, \nu \in {\mathcal P}(n),$$
where ${\mathcal P}(n)$ is the set of partitions of $n$.
Let
$$\mathsf{Z}'_{GW}(\Pp\times \C^2)_{n[\Pp],\lambda,\mu,\nu}\in {\mathbb Q}(t_1,t_2)((u))$$
be the reduced Gromov-Witten
partition function, see Section 3.2 of \cite{jbrp}.

The $T$-equivariant cohomology of $\Hilb_n$ has
a canonical Nakajima basis indexed by ${\mathcal P}(n)$.
Define the series $\langle \lambda,\mu,\nu\rangle^{\Hilb_n}$ of
 3-point invariants by a sum over curve degrees:
$$\langle \lambda,\mu,\nu
\rangle^{\Hilb_n} = \sum_{d\geq 0} q^d
\langle \lambda,\mu,\nu
\rangle^{\Hilb_n}_{0,3,d}.$$

The results of \cite{jbrp} together with 
our calculation of the $3$-point functions of the
Hilbert scheme yields a {\em Gromov-Witten/Hilbert
correspondence} discussed in Section \ref{s_relGW}.

\vspace{+10pt}
\noindent{\bf Theorem.} After the variable change $e^{iu}=-q$,
$$(-iu)^{-n+\ell(\lambda)+\ell(\mu)+\ell(\nu)} \mathsf{Z}'_{GW}(\Pp\times \C^2)_{n[\Pp],\lambda,\mu,\nu} 
= (-1)^n\langle \lambda,\mu,\nu
\rangle^{\Hilb_n}.$$
\vspace{+5pt}

\noindent In fact, our Hilbert scheme calculations were motivated
by the correspondence.

The Gromov-Witten and Donaldson-Thomas theories of
$\Pp \times \C^2$ are related by the correspondence
conjectured in \cite{mnop1,mnop2} and refined for the $T$-equivariant
context in \cite{jbrp}. 
The quantum cohomology of $\Hilb_n$ provides
a third vertex of equivalence:

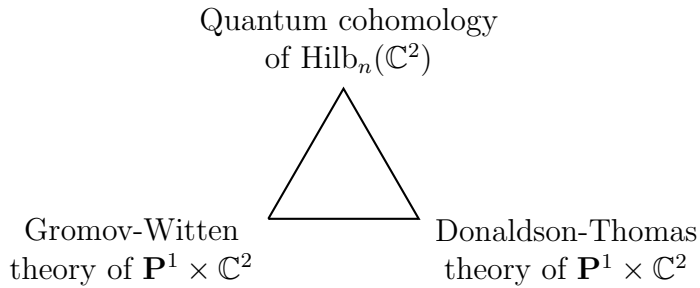
\begin{figure}[hbtp]\psset{unit=0.5 cm}
  \begin{center}
    \begin{pspicture}(-6,-2)(10,6)
    \psline(0,0)(2,3.464)(4,0)(0,0) 
    \rput[rt](0,0){
        \begin{minipage}[t]{3.64 cm}
          \begin{center}
            Gromov-Witten \\ theory of $\Pp \times \C^2$
          \end{center}
        \end{minipage}}
    \rput[lt](4,0){
        \begin{minipage}[t]{3.64 cm}
          \begin{center}
             Donaldson-Thomas\\ theory of $\Pp \times \C^2$
          \end{center}
        \end{minipage}}
    \rput[cb](2,4.7){
        \begin{minipage}[t]{4 cm}
          \begin{center}
           Quantum cohomology \\ of $\Hilb_n(\C^2)$
          \end{center}
        \end{minipage}}
    \end{pspicture}
  \end{center}
\end{figure}

\noindent The proof of the triangle of equivalences is completed in
 \cite{GWDT}.{\footnote{The equivalences of the triangle
have been extended to the $A_n$ case in \cite{M,MO1,MO2}.}}

The role played by the quantum differential equation \eqref{qode} 
in Gromov-Witten and Donaldson-Thomas theories is the following. 
The fundamental solution of \eqref{qode} is related, on the 
one hand, to general triple Hodge integrals on the moduli 
space of curves and, on the other hand, to 
the equivariant vertex \cite{mnop1,mnop2}
with one infinite leg, see \cite{GWDT}.

\subsection{Acknowledgments} 
We thank J. Bryan, B. Fantechi, T. Graber, N. Katz, J. Koll\'ar, M. Manetti,
D. Maulik,  N. Nekrasov, and Z. Ran for many  
valuable discussions. Both authors were partially 
supported by the Packard 
foundation and the NSF.

\section{The operator $\MM_D$}\label{s2} 

\subsection{Fock space formalism}
\label{fs}

We review the Fock space description of the $T$-equivariant 
cohomology of the Hilbert scheme of points of $\C^2$, see \cite{Gron, Nak}. By definition, the Fock space $\cF$ 
is freely generated over $\Q$ by commuting 
creation operators $\alpha_{-k}$, $k\in\Z_{>0}$,
acting on the vacuum vector $\vac$. The annihilation 
operators $\alpha_{k}$, $k\in\Z_{>0}$, kill the vacuum 
$$
\alpha_k \cdot \vac =0,\quad k>0 \,,
$$
and satisfy the commutation relations
$$
\left[\alpha_k,\alpha_l\right] = k \, \delta_{k+l}\,. 
$$

A natural basis of $\cF$ is given by 
the vectors  
\begin{equation}
  \label{basis}
  \lv \mu \rang = \frac{1}{\zz(\mu)} \, \prod \alpha_{-\mu_i} \, \vac \,.
\end{equation}
indexed by partitions 
$\mu$. Here, $$\zz(\mu)=|\Aut(\mu)| \, \prod \mu_i$$ is the usual 
normalization factor. 
Let the length $\ell(\mu)$ denote the number of 
parts of the partition $\mu$.

The {\em Nakajima basis} defines a canonical isomorphism,
\begin{equation}
\cF \otimes _{\mathbb Q} {\mathbb Q}[t_1,t_2]\stackrel{\sim}{=} 
\bigoplus_{n\geq 0} H_{T}^*(\Hilb_n,{\mathbb Q}).
\label{FockHilb}
\end{equation}
The Nakajima basis element corresponding to  
$\lv \mu \rang$  is
$$\frac{1}{\Pi_i \mu_i} [V_\mu]$$
where $[V_\mu]$ is (the cohomological dual of)
the class of the subvariety of $\Hilb(\C^2,|\mu|)$
with generic element given by a union of 
schemes of lengths $$\mu_1, \ldots, \mu_{\ell(\mu)}$$ supported
at $\ell(\mu)$ distinct points of $\C^2$. 
The vacuum vector $\vac$ corresponds to the unit in 
$H_T^*(\Hilb_0, {\mathbb Q}).$
As before, $t_1,t_2$ are the equivariant parameters corresponding
to the weights of the $T$-action on $\C^2$.

The subspace of $\cF\otimes_{\mathbb Q} {\mathbb Q}[t_1,t_2]$ corresponding to $H^*_T(\Hilb_n,{\mathbb Q})$
is spanned by the vectors \eqref{basis} with $|\mu|=n$. The subspace can also 
be described as the $n$-eigenspace of the {\em  energy operator}: 
$$
|\cdot| = \sum_{k>0} \alpha_{-k} \, \alpha_k \,.
$$
The vector $\lv 1^n \rang$ corresponds by to the identity in $H^*_T(\Hilb_n,{\mathbb Q})$.
A straightforward calculation shows 
$$
D = - \lv 2,1^{n-2} \rang \,.
$$

The standard inner product on the $T$-cohomology  induces the following 
{\em nonstandard} inner product on Fock space after an extension of scalars:
\begin{equation}
  \label{inner_prod}
  \lang \mu | \nu \rang = 
\frac{(-1)^{|\mu|-\ell(\mu)}}{(t_1 t_2)^{\ell(\mu)}} 
\frac{\delta_{\mu\nu}}{\zz(\mu)} \,. 
\end{equation}
With respect to the inner product, 
\begin{equation}
  \label{adjoint}
  \left(\alpha_{k}\right)^* = (-1)^{k-1} (t_1 t_2)^{\sgn(k)} \, 
\alpha_{-k}\,.
\end{equation}

\subsection{Main theorem}
\label{t1}

The following operator 
on Fock space plays a central role in the paper: 
\begin{multline}
  \label{theM} 
\MM(q,t_1,t_2) = (t_1+t_2) \sum_{k>0} \frac{k}{2} \frac{(-q)^k+1}{(-q)^k-1} \,
 \alpha_{-k} \, \alpha_k  + \\
\frac12 \sum_{k,l>0} 
\Big[t_1 t_2 \, \alpha_{k+l} \, \alpha_{-k} \, \alpha_{-l} -
 \alpha_{-k-l}\,  \alpha_{k} \, \alpha_{l} \Big] \,.
\end{multline}
The $q$-dependence of $\MM$ is only in the first sum in \eqref{theM}
which acts diagonally in the basis \eqref{basis}. The two terms in the 
second sum in \eqref{theM} are known respectively as the splitting and joining terms. 
The operator $\MM$ commutes with the energy operator 
$|\cdot|$, and 
\begin{equation}
  \label{Madj}
  \MM^* = \MM 
\end{equation}
with respect to \eqref{adjoint}. 

\begin{Theorem}\label{T1}
Under the identification \eqref{FockHilb}, 
  \begin{equation}
    \label{MM_D}
    \MM_D = \MM - \frac{t_1+t_2}{2} \, \frac{(-q)+1}{(-q)-1} \, 
|\,\cdot\,| \,\,
  \end{equation}
is the operator of 
small quantum multiplication by the divisor $D$ in the
$T$-equivariant cohomology of the Hilbert scheme of points of $\C^2$. 
\end{Theorem}

\begin{Corollary}\label{cor1} 
The divisor class $D$ generates the small quantum ring
$$QH^*_T(\Hilb_n,{\mathbb Q})$$ 
over $\Q(q,t_1,t_2)$.   
\end{Corollary}

In the basis \eqref{basis}, 
the matrix elements of $\MM_D$ are {\em integral} ---  
the matrix elements lie in $\Z[t_1,t_2][[q]]$. 


\subsection{Calogero-Sutherland operator}

The classical
multiplication by the divisor $D$ and the connection to the Calogero-Sutherland
operator, 
\begin{equation}
  \label{HCS}
  \bH_{CS}= \frac12 \sum_i 
\left( z_i \frac{\partial}{\partial z_i}\right)^2 +
\theta(\theta-1) \sum_{i<j} \frac1{|z_i-z_j|^2} \,,
\end{equation}
are recovered by setting $q=0$ in $\MM_D$.

The 
operator $\bH_{CS}$ describes quantum-mechanical particles on the torus 
$|z_i|=1$ interacting via the potentials $|z_i-z_j|^{-2}$. The parameter
$\theta$ adjusts the strength of the interaction. 
The function 
$$
\phi(z) = \prod_{i<j} (z_i-z_j)^\theta
$$
is an eigenfunction of $\bH_{CS}$, and the operator
$\phi \, \bH_{CS} \, \phi^{-1}$ preserves the space of symmetric
polynomials in the variables $z_i$. Therefore,
via the identification 
$$
p_\mu(z) = \zz(\mu) \, \lv \mu \rang \,,
$$
where
$$
p_\mu(z) = \prod_k \sum_i z_i^{\mu_k} \,,
$$
the operator
$\phi \, \bH_{CS} \, \phi^{-1}$
acts on  Fock space.

A direct computation shows  the operator  $\phi \, \bH_{CS} \, \phi^{-1}$ equals
\begin{equation}
 \Delta_{CS}= \frac{1-\theta}{2}\sum_k k \, \alpha_{-k} \, \alpha_k  + \\ 
\frac12 \sum_{k,l>0} 
\Big[\alpha_{-k-l} \, \alpha_{k} \, \alpha_{l} + \theta \, 
 \alpha_{k+l}\,  \alpha_{-k} \, \alpha_{-l} \Big] \, 
\end{equation}
modulo scalars and a multiple of the momentum 
operator $\sum_i z_i \frac{\partial}{\partial z_i}$, see \cite{St}.
We find
\begin{equation}
\label{MDCS}
  \MM(0) = 
- \,t_1^{\ell(\,\cdot\,)+1} \, 
\Delta_{CS}\Big|_{\theta=-t_2/t_1} \, t_1^{-\ell(\,\cdot\,)}  \ .
\end{equation}
The well-known duality $\theta\mapsto1/\theta$
in the Calogero-Sutherland model corresponds to 
the permutation of $t_1$ and $t_2$.

\subsection{Eigenvectors}\label{s_eigM} 
Let $\lambda$ be a partition of $n$. Let $\cI_\lambda$
denote the associated $T$-fixed (monomial) ideal,
\begin{equation}
  \label{cJla}
 \cI_\lambda=\big( x^{j-1} y^{i-1} \big)_{\square=(i,j)\notin\lambda} \subset 
\C[x,y] \,.
\end{equation}
The map
$\lambda \mapsto \cI_\lambda$
is a bijection between the set of  partitions ${\mathcal P}(n)$  and
the set of 
$T$-fixed points 
$\Hilb_n^T \subset \Hilb_n$.

The eigenvectors of the classical multiplication by $D$ in $H^*_T(\Hilb_n,{\mathbb Q})$
are the classes of the $T$-fixed points of $\Hilb_n$,
$$
\left[\cI_\lambda\right] \in H^{2n}_T(\Hilb_n, {\mathbb Q})\,, \quad \lambda\in {\mathcal P}(n)\,.
$$
The eigenvalues of $\MM_D(0)$ are determined by the functions
\begin{equation}
  \label{def_c}
 - c(\lambda;t_1,t_2) = - \sum_{(i,j)\in \lambda} \Big[
(j-1) t_1 + (i-1) t_2 \Big] \,,
\end{equation}
for $\lambda \in {\mathcal{P}}(n)$. 
The sum in \eqref{def_c} is the trace of the $T$-action on
 $\C[x,y]/\cI_\lambda$. 

We will denote by 
$$
\bJ^\lambda\in \cF\otimes \Q[t_1,t_2] 
$$
the image of $\left[\cI_\lambda\right]$ in Fock space. The corresponding symmetric 
function is the {\em integral form} of the Jack polynomial \cite{Mac}.

\section{Proof of Theorem \ref{T1}}

\subsection{3-point functions}

We must prove the matrix elements of $\MM_D$ yield the $T$-equivariant
3-point functions of the Hilbert scheme:
\begin{equation}
  \label{tMM_Hilb}
  \lang \mu \Bv \tMM \Bv \nu \rang = 
\sum_{d\ge 0} q^d \, 
\blang \mu, D , \nu \brang^{\Hilb_n}_{0,3,d} \,,
\end{equation}
where $\mu,\nu \in {\mathcal P}(n)$.
The matrix elements on the left side of \eqref{tMM_Hilb} are calculated with the nonstandard
inner product \eqref{inner_prod}.

The $q^d$ coefficients of
the left and right sides of \eqref{tMM_Hilb} will be denoted by the respective brackets:
$$\lang \mu \Bv \tMM \Bv \nu \rang_d, \ \  \lang \, \mu, D, \nu \rang_d.$$ 
The calculation of the $T$-equivariant (classical) cohomology of $\Hilb_n$
implies the equality \eqref{tMM_Hilb} in degree 0.

\subsection{Definitions}

Though  $\Hilb_n$ is not compact, the $T$-equivariant 
Gromov-Witten invariants are well-defined.
The $T$-fixed locus of the moduli space of maps to $\Hilb_n$ is
a proper Deligne-Mumford stack. The $T$-equivariant Gromov-Witten theory
may be defined by a residue integral on the $T$-fixed locus via the
virtual localization formula \cite{GrabPan}.

An equivalent geometric definition of the $T$-equivariant
Gromov-Witten theory of $\Hilb_n$ is obtained in the {fixed point} basis.
The classes of the $T$-fixed points of $\Hilb_n$ span a basis of
the (localized) equivariant cohomology of $\Hilb_n$.
The locus of maps to $\Hilb_n$ meeting a $T$-fixed point {\em is} compact.
For example, the 3-point functions in the $T$-fixed point basis are:
$$\lang [\cI_\mu], [\cI_\nu], [\cI_\xi] \rang_{d} = 
\int_{[\overline{M}_{0,3}(\Hilb_n,d)]^{vir}} 
\ \text{ev}_1^*([\cI_\mu]) \cup \text{ev}_2^*([\cI_\nu]) \cup \text{ev}_3^*([\cI_\xi]),$$
where the integral sign denotes $T$-equivariant push-forward to a point.

By either definition, the $T$-equivariant Gromov-Witten invariants of
$\Hilb_n$ with insertions in $H^*_{T}(\Hilb_n,{\mathbb Q})$
have values in the ring $\Q(t_1,t_2)$.

\subsection{The Nakajima basis revisited}
The Nakajima basis was defined in Section \ref{fs} with respect to the
identity element of $H^*_T(\C^2,{\mathbb Q})$.
We may also define a Nakajima basis with respect to the class of
the origin
$$
[\mathbf{0}] = t_1 t_2 \in H^{4}_T(\C^2, {\mathbb Q}).  
$$
The Nakajima basis at the origin is determined by:
\begin{equation}\label{ffxx}
\lv \mu([\mathbf{0}]) \rang= (t_1 t_2)^{\ell(\mu)} \lv \mu \rang
\in H^{2(|\mu|+\ell(\mu))}_T(\Hilb_{|\mu|}, {\mathbb Q})\,. 
\end{equation}
By linearity,
\begin{align}
\blang \mu, D, \nu \brang_{d} &= 
\frac1{(t_1 t_2)^{\ell(\mu)}} \, 
\blang \mu([\mathbf{0}]), D, \nu \brang_{d} 
 \label{C^2_P^2}
\\
&= 
\frac1{(t_1 t_2)^{\ell(\nu)}} \, 
\blang \mu, D, \nu ([\mathbf{0}]) \brang_{d} \,. 
\notag
\end{align}

The Gromov-Witten invariants on the right side of
\eqref{C^2_P^2} are intersection products in a compact space (because
of the occurrences of the class $[{\mathbf 0}]$) and,
therefore, take values in $\Q[t_1,t_2]$. In order for such
Gromov-Witten invariants to be nonzero, the total codimension of
the insertions must not be less than the virtual dimension of
the moduli space of maps.

The virtual dimension of the moduli space of genus 0 maps to
$\Hilb_n$ with 3 marked points is:
\begin{eqnarray*}
\text{vir dim}_{\C} \ \overline{M}_{0,3}(\Hilb_n,\beta)&  = &
\int_\beta c_1(T_{\Hilb_n}) + 2n-3 + 3 \\
& = & 2n.
\end{eqnarray*}
Since $\Hilb_n$ is  holomorphic symplectic, the first Chern class
of the tangent bundle is trivial.

The (complex) codimensions of the insertions
$\mu([{\mathbf 0}])$, $D$, and $\nu$ sum to
$$ |\mu|+ \ell(\mu) + 1 + |\nu| -\ell(\nu) = 2n+1+ \ell(\mu)-\ell(\nu).$$
Similarly, the codimensions of the insertions
$\mu$, $D$, $\nu([{\mathbf 0}])$ sum to
$$2n+1+ \ell(\nu)-\ell(\mu).$$
The inequality
\begin{equation}
  \label{nonvineq}
  |\ell(\mu)-\ell(\nu)|\le 1
\end{equation}
is therefore a necessary condition for the nonvanishing 
of the $3$-point invariants \eqref{tMM_Hilb}.

If $\ell(\mu)\ne \ell(\nu)$, then one of the invariants
on the right in \eqref{C^2_P^2} has degree 0 in $t_1$ and $t_2$ and, 
therefore, is purely classical by  
Lemma \ref{Pt1t2} below. Hence, for  $\ell(\mu)\ne \ell(\nu)$,
equation \eqref{tMM_Hilb} 
reduces to the known formula for the classical 
multiplication by $D$.

\subsection{Reduced virtual classes}
\label{rvclass}
\subsubsection{A fixed domain $C$}
\label{aaa}
Let $C$ be fixed, pointed, nodal, genus $g$ curve. Let
$$M_C(\Hilb_n,d)$$ denote the moduli space of maps from $C$ to $\Hilb_n$ of
degree $d>0$. 
Let 
$$\pi: C \times M_C(\Hilb_n,d) \rightarrow M_C(\Hilb_n,d)$$
denote the projection, and let
$$f:  C \times M_C(\Hilb_n,d) \rightarrow \Hilb_n$$
denote the universal map.
The canonical morphism
\begin{equation}\label{reld}
R^\bullet \pi_*(f^* T_{\Hilb_n})^\vee \rightarrow L_{M_C}^\bullet
\end{equation}
determines a perfect obstruction theory on $M_C(\Hilb_n,d)$, see \cite{Beh,BehFan,LiTian}.
Here, $L_{M_C}^\bullet$ denotes the cotangent complex of $M_C(\Hilb_n,d)$.

Let $dx\wedge dy$ be the standard holomorphic symplectic form
on $\C^2$. Let $$\C dx\wedge dy$$ be the associated 1-dimensional $T$-representation
of weight $-(t_1+t_2)$. The form $dx\wedge dy$
induces a canonical holomorphic symplectic form $\gamma$ on
$\Hilb_n$. The $T$-representation $\C \gamma$ has weight
$-n(t_1+t_2)$.

Let $\Omega_\pi$ and $\omega_\pi$ denote respectively the sheaf
of relative differentials and the relative dualizing sheaf.
There is a canonical map
$$f^*(\Omega_{\Hilb_n}) \rightarrow \Omega_\pi \rightarrow \omega_\pi.$$
After dualizing, we obtain
\begin{equation}\label{fgt}
\omega_\pi^* \rightarrow f^*(T_{\Hilb_n}).
\end{equation}
The map \eqref{fgt} and the holomorphic symplectic form $\gamma$ on $\Hilb_n$ together yield
a map
$$f^*(T_{\Hilb_n}) \rightarrow \omega_\pi \otimes (\C\gamma)^\vee.$$
We obtain
$$R^\bullet\pi_*(\omega_\pi)^\vee \otimes \C \gamma \rightarrow R^\bullet\pi_*(f^*T_{\Hilb_n})^\vee.$$
Finally, we consider the induced cut-off map
$$\iota: \tau_{\leq -1} R^\bullet\pi_*(\omega)^\vee \otimes \C \gamma  \rightarrow
R^\bullet\pi_*(f^*T_{\Hilb_n})^\vee.$$

The complex $\tau_{\leq -1} R^\bullet\pi_*(\omega)^\vee\otimes \C\gamma$ is represented
by a trivial bundle of rank $1$ with representation $\C\gamma$ in degree $-1$.
Consider the mapping cone $C(\iota)$ of $\iota$.
Certainly $R^\bullet\pi_*(f^*T_{\Hilb_n})^\vee$ is represented
by a two term complex. An elementary argument using the positive degree $d>0$
condition shows
the complex $C(\iota)$ is also two term.

By Ran's results{\footnote{The required deformation theory can also
be found in a recent paper by M. Manetti \cite{Man}.
The comments of Sections \ref{aaa} and \ref{bbb} apply equally
to Manetti's results.}} on deformation theory and the semiregularity map, 
there is a canonical map
\begin{equation}
\label{rvc}
C(\iota) \rightarrow L_{M_C}^\bullet
\end{equation}
induced by \eqref{reld}, see \cite{zzz2}.
Ran proves the obstructions to deforming maps from $C$ to
a holomorphic symplectic manifold lie in the kernel
of the semiregularity map. After dualizing,  Ran's result
precisely shows \eqref{reld} factors through the cone $C(\iota)$.

The map \eqref{rvc} defines a {\em new} perfect obstruction theory
on $M_C(\Hilb_n,d)$. The conditions of cohomology isomorphism in degree 0 and
the cohomology surjectivity in degree $-1$ are both induced from
the perfect obstruction theory \eqref{reld}.

We view \eqref{reld} as the {\em standard} obstruction theory and
\eqref{rvc} as the {\em reduced} obstruction theory.
Both obstruction theories are $T$-equivariant since the 
morphism of complexes involved are $T$-equivariant.

\subsubsection{Ran's results}
\label{bbb}
Two aspects of the application of Ran's 
deformation results here warrant further comment. 

First, a main
technical advance in \cite{zzz2} is the study of obstructions for deformations over Artin local
rings: the case of deformations over the curvilinear schemes $\C[\epsilon]/(\epsilon^n)$
was treated earlier in \cite{bloch,zzz1}. The Artin local case is needed here.

Second, Ran's proof requires a nonsingular projective target variety with a holomorphic
symplectic form. While $\Hilb_n$ is {\em not} complete, Ran's argument can be nevertheless
be applied by the following construction. Let 
$$f:C \rightarrow \Hilb_n$$
be a stable map. The image of $f$ under composition with the
Hilbert/Chow morphism,
$$\rho_{HC}: \Hilb_n \rightarrow \text{Sym}^n \C^2,$$
must be a point $$\sum_{i=1}^{\ell(\mu)} \mu_i[p_i].$$
We may view $f$ as a map to the fiber $\rho_{HC}^{-1}(\sum \mu_i[p_i])$ of the Hilbert/Chow morphism.

Let $S$ be a nonsingular, projective, $K3$ surface, and let $q_1, \ldots, q_{\ell(\mu)}$
be distinct points.
Consider the fiber $$\rho_{HC}^{-1}\left(\sum\mu_i[q_i]\right) \subset \Hilb(S,n).$$
 Local analytic charts on $\C^2$ and $S$ at the points $p_i$ and $q_i$ induce
a local analytic isomorphism  of $\Hilb_n$ and $\Hilb(S,n)$
in a neighborhood of the two fibers of $\rho_{HC}$.
Hence, the deformation theory of $f$ over Artin local rings can be 
studied on $\Hilb(S,n)$. Since $S$ and $\Hilb(S,n)$ are holomorphic
symplectic, Ran's results imply the obstructions lie in the
kernel of the semiregularity map for $\Hilb(S,n)$. The statement
implies precisely the required deformation theory statement for $\Hilb_n$.

 \subsubsection{The reduced absolute theory}
The results of Sections \ref{aaa}-\ref{bbb} define a $T$-equivariant reduced
obstruction theory of maps to $\Hilb_n$ relative to the Artin stack
${\mathfrak M}$ of pointed genus $g$ curves. A $T$-equivariant reduced absolute theory is
obtained via a distinguished triangle in the usual way, see \cite{Beh,BehFan,LiTian}.

Since the new obstruction theory differs from the standard theory by the
1-dimensional obstruction space $(\C \gamma)^\vee$, we find
\begin{eqnarray*}
[\overline{M}_{g,n}(\Hilb_n,d)]^{vir}_{s}& =& 
c_1\left(  (\C \gamma)^\vee \right) \cap [\overline{M}_{g,n}(\Hilb_n,d)]^{vir}_{r}  \\
& = & (t_1+t_2) \cdot [\overline{M}_{g,n}(\Hilb_n,d)]^{vir}_{r},
\end{eqnarray*}
for $d>0$.
Here, $s$ and $r$ denote the standard and reduced theories.

\begin{Lemma}\label{Pt1t2}
$T$-equivariant Gromov-Witten invariants of $\Hilb_n$ of positive degree with
insertions in $H_T^*(\Hilb_n,{\mathbb Q})$ are divisible by $t_1+t_2$. 
\end{Lemma}

\begin{proof}
The $T$-equivariant Gromov-Witten invariants lie in $\Q(t_1,t_2)$. Divisibility 
is defined by  positive valuation at
 $t_1+t_2$.

Consider the $T$-equivariant Gromov-Witten theory of $\Hilb_n$ in the
Nakajima basis at the origin.
By compactness, the invariants lie in $\Q[t_1,t_2]$. By the
construction  of the reduced virtual class, the invariants are
divisible by $t_1+t_2$.

The Nakajima basis with respect to the identity spans $H_T^*(\Hilb_n,{\mathbb Q})$
as a ${\Q}$-vector space. The relation \eqref{ffxx} concludes the proof.
\end{proof}

\subsection{Additivity}

Denote the reduced
invariants of $\Hilb_n$ by curved brackets: 
$$
\blang \mu([\mathbf{0}]), D, \nu \brang_{d} =
(t_1+t_2) \,\dlang \mu([\mathbf{0}]), D, \nu \drang_{d} \,, 
\quad d>0 \,. 
$$
If $\ell(\mu)= \ell(\nu)$, the integral
$\dlang \mu([\mathbf{0}]), D, \nu \drang_{d}$
is a nonequivariant constant.
Let $$\xi=\sum \mu_i [p_i]\  \in \Sym^n \C^2,$$
where  $\{p_1,\dots,p_{\ell(\mu)}\}\subset\C^2$
are distinct points.
By the definition of the Nakajima basis,
we can replace the equivariant class $ \mu([\mathbf{0}])$ in the integrand by the 
nonequivariant class 
\begin{equation}\label{g345}
\mu(\xi)= \frac{1}{\zz(\mu)} \, 
\rho_{HC}^{-1}\left(\xi\right)
\end{equation}
where
$$
\rho_{HC}: \Hilb_n \to \Sym^n \C^2
$$
is the Hilbert/Chow morphism as before.

Every rational curve in $\Hilb_n$  is contracted by the 
Hilbert/Chow morphism. The moduli space of maps 
connecting the locus $\mu(\xi)$ and the Nakajima cycle $\nu$ is 
 isomorphic to the moduli space of stable 
maps to the product
\begin{equation}\label{ffvvs}
\prod_{i=1}^{\ell(\mu)} \Hilb_{|\mu_i|,p_i}
\end{equation}
in case $\mu=\nu$ and empty otherwise.
Here, $\Hilb_{m,p}\subset \Hilb_m$ denotes the subspace of schemes supported at $p$. 

The moduli space of maps to the product,
\begin{equation}\label{dr555}
\overline{M}_{0,3}\left(\prod_{i=1}^{\ell(\mu)} \Hilb_{|\mu_i|,p_i}, d\right),
\end{equation}
has components corresponding to the
different distributions of the total degree $d$ 
among the factors. Consider a component 
$$\overline{M}[j,k] \subset \overline{M}_{0,3}\left(\prod_{i=1}^{\ell(\mu)} \Hilb_{|\mu_i|,p_i}, d\right),$$
for which the
degree splitting has at least two non-zero terms corresponding to
the points $p_j$ and $p_k$.

The moduli space $\overline{M}[j,k]$ has a standard obstruction theory obtained
from the standard obstruction theory of $\overline{M}_{0,3}(\Hilb_n, d)$.
The standard obstruction theory of $\overline{M}[j,k]$
has a 2-dimensional quotient 
obtained from the 2-dimensional family of holomorphic symplectic
forms (at $p_j$ and $p_k$). Exactly following the construction
of Section \ref{rvclass}, we obtain a {\em doubly} reduced obstruction theory 
by reducing the obstruction space by the 2-dimensional quotient.
The nonequivariant integral of the (singly) reduced theory over such a
component vanishes since the singly reduced theory contains
an additional 1-dimensional trivial factor.

We conclude the only components of \eqref{dr555}
which contribute to the integral
$$\dlang \mu(\xi), D, \mu \drang_{d}$$
are those for which the degree $d$ is distributed entirely to a 
single factor of the product. 

After unraveling the
definitions, we obtain three basic results governing the
insertion $D$:
\begin{enumerate}
\item[(i)] 
$\blang \mu, D, \nu \brang_{d>0} =0$
for $\mu \neq \nu$,
\item[(ii)] 
$\blang \mu, D, \mu \brang_d = 
{\gamma_{n,d}}{(t_1t_2)^{-\ell(\mu)}}(t_1+t_2)$ 
where $\gamma_{n,d}\in {\mathbb Q}$,
\item[(iii)]
the {\em addition
formula},
\begin{equation}
  \label{addm}
\frac{\blang \mu, D, \mu \brang_{d>0}}
{\lang \mu \bv \mu \rang} =
\sum_i 
\frac{\blang \mu_i, D, \mu_i \brang_{d>0}}
{\lang \mu_i \bv \mu_i \rang} \,. 
\end{equation}
\end{enumerate}
The corresponding properties of $\MM_D$, including the
addition formula,
\begin{equation}
  \label{addmM}
\frac{\blang \mu, \MM_D, \mu \brang_{d>0}}
{\lang \mu \bv \mu \rang} =
\sum_i 
\frac{\blang \mu_i,\MM_D, \mu_i \brang_{d>0}}
{\lang \mu_i \bv \mu_i \rang} \,, 
\end{equation}
are directly verified.

\subsection{Induction strategy}

We will prove Theorem 1 by induction on $n$. If $n=0,1$, the
operator $\MM_D$ vanishes. The insertion $D$ is 0 for $n=0,1$,
so Theorem 1 is valid.

Let $n>1$. We proceed by induction on the degree $d$.
The induction step relies upon the  
addition formulas \eqref{addm}-\eqref{addmM}.
For each degree $d\geq1$, we will compute a 
 3-point invariant
$$
\blang \gamma_1, D , \gamma_2 \brang_{d}$$
for which the expansions of the classes
$$\gamma_1,\gamma_2 \in H^{4n}_T(\Hilb_n,{\mathbb Q}),$$
in the Nakajima basis contain nontrivial multiples 
{\em not divisible by $(t_1+t_2)$} of the class $|n\rangle$.
By the addition rules, if 
\begin{equation}
  \label{t123}
  \lang \gamma_1 \Bv \tMM \Bv \gamma_2 \rang_d = 
\blang \gamma_1, D , \gamma_2 \brang_{d} \,,
\end{equation}
then \eqref{tMM_Hilb} is proven for $\Hilb_n$ in degree $d$. 

Both sides of 
\eqref{t123} are constant multiples of $t^{2n}_1(t_1+t_2)$ modulo
$(t_1+t_2)^2$. Since $\blang (n), D, (n) \brang_d$ is determined
by the constant $\gamma_{n,d}$ where
$$\blang (n), D, (n) \brang_d = -\frac{\gamma_{n,d}}{t_1^2} (t_1+t_2)\ \mod
(t_1+t_2)^2,$$
we need only  
verify the equality (\ref{t123}) modulo $(t_1+t_2)^2$.

\subsection{Induction step: I}
Let $n>1$ and let $d\geq 1$.
For the induction step, we will compute the invariant
\begin{equation}
  \label{inv_comp}
  \lang \left[\cI_{(n)}\right], D, \left[\cI_{(n-1,1)}\right]
\rang_d \,.
\end{equation}
Following the notation of Section \ref{s_eigM},
 $\cI_{\lambda}$ denotes the monomial ideal corresponding to the
partition $\lambda$, and $[\cI_\lambda]$ denotes the
$T$-equivariant class of the associated fixed point in $\Hilb_{|\lambda|}$.

The $T$-fixed point $[\cI_{\lambda}]$
corresponds to the Jack polynomial $$\bJ^\lambda\in\cF\otimes \Q[t_1,t_2].$$ 
For $\theta=-t_2/t_1=1$, the Jack polynomials specialize to the Schur functions.
Hence,
$$
\bJ^\lambda \equiv \frac{(-1)^{|\lambda|} \, |\lambda|!}{\dim \lambda}
\sum_\mu \chi^\lambda_\mu  \, t_1^{|\lambda|+\ell(\mu)}\, 
\bv \mu \brang  \mod t_1+t_2 \,,
$$
where $\dim \lambda$ is the dimension of the
representation $\lambda$ of the symmetric group and 
$\chi^\lambda_\mu$ is the associated character evaluated
on the conjugacy class $\mu$, see \cite{Mac}. 
In particular, the coefficient of $\bv n \brang$
in the expansion of both $\bJ^{(n)}$ and $\bJ^{(n-1,1)}$
is nonzero. 

The operator $\tMM-\tMM(0)$ formed by the 
terms of positive $q$ degree  in the operator $\tMM$ 
acts diagonally in the basis $\bv \mu \brang$. 
Since $\chi^{(n)}$ is the trivial character and
$$
\dim (n-1,1)=n-1\,,
$$ 
we conclude 
\begin{multline*}
\lang \bJ^{(n)} \Bv \tMM-\tMM(0) \Bv  \bJ^{(n-1,1)} \rang
\equiv \\
(-1)^n (t_1+t_2) \,\frac{t_1^{2n} \, (n!)^2}{n-1} 
\left(\chi^{(n-1,1)},F\right)_{L^2(S(n))}
 \mod  (t_1+t_2)^2\,,
\end{multline*}
where $F$ is the function on the symmetric group $S(n)$ taking the  
value
\begin{equation}
F(\mu) = -|\mu| \frac{q}{1+q} - \sum_i \mu_i^2 \frac{(-q)^{\mu_i}}{1-(-q)^{\mu_i}}
\label{Fmu}
\end{equation}
on a permutation with cycle type $\mu$, and $(,)$
is the standard inner product on $L^2(S(n))$ with respect to 
which the characters are orthonormal. 

The first term in \eqref{Fmu} is a constant function and, hence, 
orthogonal to any nontrivial character. The second term can be written
as
$$
\sum_i \mu_i^2 \frac{(-q)^{\mu_i}}{1-(-q)^{\mu_i}}= \sum_{k\ge 1} f''(\mu,(-q)^k)\,,
$$
where
$$
f(\mu,z) = \sum_i z^{\mu_i}$$
and differentiation is taken with respect to the operator $z\frac{d}{dz}$,
$$
 \quad f''(\mu,z)=\left(z\frac{d}{dz}\right)^2 f(\mu,z) \,.
$$
The evaluation of the inner product $ \left(\chi^{(n-1,1)},f''\right)$
is obtained by differentiating the
$$
(a,b,c)=(n-1,1,0)
$$
case of the following result.

\begin{Lemma}\label{aapo}
The Fourier coefficients of  $f'$ are: 
  \begin{equation}
\label{harm_dec}
    \left(\chi^{\lambda},f'\right)_{L^2(S(n))} =
    \begin{cases}
      \sum_{k=1}^{a} z^k\,, & \lambda=(a) \,,\\
     (-1)^{c+1} \,z^{a+c+1} + (-1)^{c} \,z^{b+c} \,,
& \lambda=(a,b,1^c)\,,\\
0\,, & \textup{otherwise}\,. 
    \end{cases}
  \end{equation}
\end{Lemma}

\begin{proof}
The inner product $\left(\chi^{\lambda},f\right)$ is the image of 
the Schur function $s_\lambda$ under the map on symmetric functions
induced by following transformations of the power-sums:
$$p_\mu \mapsto f(\mu,z).$$ 
We observe
$$
p_\mu(1,\underbrace{z,\dots,z}_\textup{$N$ times}\,)=1 + N \, f(\mu,z) + O(N^2) \,.
$$
Hence, the map $p_\mu\mapsto f(\mu,z)$ is the linear coefficient in $N$
 of the above expansion. By basic properties of
the Schur functions, we find
\begin{equation}
\label{sumnu}
s_\lambda(1,\underbrace{z,\dots,z}_\textup{$N$ times}\,)=
\sum_{\nu\prec\lambda} z^{|\nu|} \, s_\nu(  \underbrace{1,\dots,1}_\textup{$N$ times}\,)
=
\sum_{\nu\prec\lambda} z^{|\nu|} \, \prod_{\square\in\nu} \frac{N+c(\square)}{h(\square)}\,,
\end{equation}
where $\nu\prec\lambda$ means that $\nu$ and $\lambda$ interlace (or, equivalently, 
that $\lambda/\nu$ is a horizontal strip). The product is over all squares $\square$
in the diagram of $\nu$, $c(\square)$ denotes the content, 
and $h(\square)$ denotes the hooklength. The linear coefficient in $N$ 
vanishes unless $\nu$ is a hook, a diagram of the form $(a,1^c)$, 
in which case the linear coefficient equals $(-1)^c/(a+c)$.

The summation over $\nu$ in \eqref{sumnu} 
telescopes to second case in \eqref{harm_dec} provided 
$\lambda$ has more than one row.  The last case in \eqref{harm_dec}
corresponds to diagrams not interlaced by a hook. 
\end{proof}

\noindent
After applying Lemma \ref{aapo}, we conclude 
\begin{multline}
\label{JJ}
\lang \bJ^{(n)} \Bv \tMM-\tMM(0) \Bv  \bJ^{(n-1,1)} \rang
\equiv \\
(-1)^n (t_1+t_2) \,\frac{t_1^{2n} \, (n!)^2}{n-1} 
\left(\frac{q}{1+q}+n\frac{(-q)^{n}}{1-(-q)^{n}}\right)
 \mod  (t_1+t_2)^2\,. 
\end{multline}
%

\subsection{Localization}
\subsubsection{Overview}
Our goal now is to reproduce the answer \eqref{JJ} by calculating the 
 $3$-point invariant
\begin{equation*}
  \lang \left[\cI_{(n)}\right], D, \left[\cI_{(n-1,1)}\right]
\rang_d \,
\end{equation*}
via localization on $\Hilb_n$.

Since the $T$-fixed locus of the moduli space $\overline{M}_{0,3}(\Hilb_n,d)$ is proper, the
 virtual localization formula of \cite{GrabPan} may be applied. However, since $\Hilb_n$
contains positive dimensional families of $T$-invariant curves, a straightforward
application is difficult. 

Our strategy for  computing  the 3-point invariant uses vanishings deduced from the
existence of the reduced obstruction theory. Since
\begin{equation}\label{cvv5} 
\lang \left[\cI_{(n)}\right], D, 
\left[\cI_{(n-1,1)}\right]\rang_d = (t_1+t_2) \dlang \left[\cI_{(n)}\right], D, \left[\cI_{(n-1,1)}\right]
\drang_d,
\end{equation}
calculation of the reduced $T$-equivariant integral on the right suffices.

Let $\TT\subset T$ denote the 1-dimensional anti-diagonal torus determined by the embedding
$$\TT\ni \xi \mapsto (\xi, \xi^{-1}) \in T.$$ Let $t$ denote the equivariant $\TT$-weight
determined by restriction,
$$t=t_1|_{\TT}= -t_2|_{\TT}.$$
Since we need to evaluate the $T$-equivariant integral
\eqref{cvv5} modulo
$(t_1+t_2)^2$, calculatation of the $\TT$-equivariant integral
$$
\dlang \left[\cI_{(n)}\right], D, \left[\cI_{(n-1,1)}\right]
\drang_d$$
suffices.

The $\TT$-fixed points of $\Hilb_n$ coincide with
the $T$-fixed points:  monomial ideals 
indexed by partitions ${\mathcal P}(n)$.
However, the $\TT$-fixed point set of the moduli space 
$\overline{M}_{0,3}(\Hilb_n,d)$ is much larger than
the $T$-fixed point set.

\subsubsection{Broken maps}
Consider the moduli space $\overline{M}_{0,k}(\Hilb_n,d)$ for $d>0$.
Let 
$$\Big[ f: C \rightarrow \Hilb_n \Big] \in \overline{M}_{0,k}(\Hilb_n,d)$$
be a $\TT$-fixed map.
If $p$ is a marking of $C$, a (fractional) $\TT$-weight $w_{p}$ is defined by the 
$\TT$-representation of the
tangent space to $C$ at $p$.
Let $P \subset C$ be a component incident to a node $s$ of $C$.
A (fractional) $\TT$-weight $w_{P,s}$ is defined by the $\TT$-representation of the
tangent space to $P$ at $s$.

We define a $\TT$-fixed map $f$
to be  {\em broken} if { either} of the following two conditions hold:
\begin{enumerate}
\item[(i)] $C$ contains a connected, $f$-contracted subcurve  $\tilde{C}$ for which the  disconnected
            curve $C\setminus \tilde{C}$ has at least two connected components which
            have positive degree under $f$.

\item[(ii)] Two non $f$-contracted components $P_1,P_2\subset C$ meet at a node $s$ of $C$ and
             have tangent weights $w_{P_1,s}$ and $w_{P_2,s}$ satisfying $w_{P_1,s}+w_{P_2,s}\neq 0$. 
\end{enumerate}
A maximal connected, $f$-contracted subcurve satisfying (i) is called a {\em breaking subcurve}.
A node satisfying (ii) is called a {\em breaking node}.
A $\TT$-fixed map which is not broken is {\em unbroken}.

A connected component of the $\TT$-fixed locus of $\overline{M}_{0,n}(\Hilb_n,d)$ is of {\em broken
type} if all the corresponding maps are broken.
Similarly, a connected component of the $\TT$-fixed locus is of {\em unbroken type} if
all the corresponding maps are unbroken.
By elementary deformation theory, every connected component is either of broken or unbroken type.

Let $\mu,\nu\in{\mathcal{P}}(n)$ be two partitions of $n$. 
A $\TT$-fixed map $[f]\in \overline{M}_{0,2}(\Hilb_n,d)$ with markings $p_1,p_2$ is said to
{\em connect} $\cI_\mu$ and $\cI_\nu$ if
$$f(p_1)= \cI_{\mu}, \ f(p_2)= \cI_{\nu},$$
and the markings do not lie on $f$-contracted components.
By definition, a connecting map must have positive degree.

\begin{Lemma}\label{444} If $f$ is  an unbroken $\TT$-fixed map of degree $d$ connecting the fixed points 
$I_\mu$ and $I_\nu$, then
$$w_{p_1} = \frac{-c(\mu;t,-t)+c(\nu;t,-t)}{d}$$
\end{Lemma}

\begin{proof} By our definitions and the stability condition, if  $f$ is unbroken then
the domain $C$ must be a chain of $r$ rational curves 
$$C= P_1 \cup \ldots \cup P_r$$ 
satisfying the condition
$$ w_{P_i,s_i}+w_{P_{i+1},s_i}=0$$
at the $i^{th}$ node $s_i$. 
Since the tangent $\TT$-representations at the fixed points of each $P_i$ have opposite weights, 
we conclude $w_{p_1}=-w_{p_2}$.

A localization calculation of the degree of the map $f$ then proves the Lemma:
\begin{eqnarray*}
d& = & \int_{f_*[C]} D \\
& = & \int_C c_1(\cO/\cI  ) \\
& = &      \frac{-c(\mu;t,-t)}{w_{p_1}} +   \sum_{i=1}^{r-1} \Big( 
\frac{ -c(f(s_i);t,-t)} {w_{P_i,s_i}}
+ \frac{ -c(f(s_i);t,-t)} {w_{P_{i+1},s_i}} \Big)
+\  \frac{-c(\nu;t,-t)}{w_{p_2}} \\
& = & \frac{-c(\mu;t,-t)}{w_{p_1}} + \frac{c(\nu;t,-t)}{w_{p_1}}.
\end{eqnarray*}
The trace of $\TT$-action on $\cO/\cI_\gamma$ 
is the function $-c(\gamma;t,-t)$, see Section \ref{s_eigM}.
\end{proof}

Let $w^d_{\mu,\nu}$ denote the tangent weight 
specified by Lemma \ref{444} at $p_1$  of an unbroken,  $\TT$-fixed, degree $d$
map connecting $\cI_\mu$ to $\cI_\nu$.
Then, $$w^d_{\nu,\mu}=-w^d_{\mu,\nu}.$$
The tangent weight $w^d_{\mu,\nu}$ is proportional to a tangent weight of $\Hilb_n$
at the fixed point $\mu$. Since the $\TT$-weights of tangent representation of $\Hilb_n$ at the
fixed points are
{\em never} $0$, we conclude the following result.

\begin{Lemma} There are no unbroken
maps connecting $\cI_\mu$ to $\cI_\mu$. 
\end{Lemma}

\subsubsection{Localization contributions}\label{slcs}

We study here 2-pointed, $\TT$-equivariant invariants of $\Hilb_n$ in positive degree, 
\begin{equation*}
\dlang [\cI_\mu], [\cI_\nu] \drang_{d} = 
\int_{[\overline{M}_{0,2}(\Hilb_n,d)]_r^{vir}} 
\text{ev}_1^*([\cI_\mu]) \cup \text{ev}_2^*([\cI_\nu]), 
\end{equation*}
where $\mu\neq \nu$.
Since the evaluation conditions lead to a proper moduli space,
the $\TT$-equivariant push-forward lies in $\Q[t]$. 


The virtual localization formula yields a sum over
the connected components of the $\TT$-fixed loci of the moduli space of maps.
Let
\begin{eqnarray*}
\dlang  [\cI_\mu], [\cI_\nu]\drang_{d}& =& 
\dlang  [\cI_\mu],[\cI_\nu]  \drang^{broken}_{d} + \dlang  [\cI_\mu], [\cI_\nu]  \drang^{unbroken}_{d},
\end{eqnarray*}
denote the separate contributions of the components of broken and unbroken type.



We index the $\TT$-equivariant localization contributions to the 2-point invariant 
$$\dlang  [\cI_\mu], [\cI_\nu]\drang_{d}$$
by graphs.
A 2-pointed tree of degree $d$ is a graph $\Gamma=(V,v_1,v_2,\rho,E,\delta)$, 
\begin{enumerate}
\item[(i)] $V$ is a finite vertex set with distinguished elements $v_1 \neq v_2$,
\item[(ii)] $\rho: V \rightarrow \Hilb_n^{\TT}$, 
\item[(iii)] $E$ is a finite edge set,
\item[(iv)] $\delta:E \rightarrow {\mathbb Z}_{>0}$ is a degree assignment,
\end{enumerate}
satisfying the following conditions
\begin{enumerate}
\item[(a)] $\Gamma$ is a connected tree,
\item[(b)] $\rho(v_1)=\cI_\mu,$ $\rho(v_2)=\cI_\nu$,
\item[(c)] if $v',v'' \in V$
are connected by an edge, then $\rho(v')\neq \rho(v'')$,
\item[(d)] if $v\neq v_1,v_2$ has edge valence 2 with neighbors $v',v''$, then
$$w^{\delta(e(v,v'))}_{\rho(v),\rho(v')} + w^{\delta(e(v,v''))}_{\rho(v),\rho(v'')}\neq 0,$$
\item[(e)] $\sum_{e\in E} \delta(e) =d.$
\end{enumerate}

Let $[f] \in \overline{M}_{0,2}(\Hilb_n,d)$ be a  $\TT$-fixed map.
We associate a 2-pointed tree, 
$$\Gamma_f=(V,v_1,v_2,\rho,E,\delta),$$ of degree $d$ to $f$ by the following
construction.
The vertex set $V$ is determined by the connected components of $f^{-1}(\Hilb_n^{\TT})$ {\em excluding} the
non-breaking nodes. In fact,
$$V=V_1\cup V_2\cup V_3,$$ is a union of three disjoint subsets:
\begin{enumerate}
\item[(1)] $V_1$ is the set of breaking subcurves,
\item[(2)] $V_2$ is the set of breaking nodes,
\item[(3)] $V_3$ is the set of nonsingular points of $C$ lying on non $f$-contracted components
           mapped to $\Hilb_n^{\TT}$.
\end{enumerate}
The two markings of $C$ are associated to distinct elements of $V_1\cup V_3$ ---
the markings  determine $v_1$ and $v_2$. The function $\rho$ is obtained from $f$.
Chains of non $f$-contracted curves of $C$ link the vertices of $V$.
The edge set $E$  is determined by such chains. The restriction of $f$ to such a chain is an unbroken connecting map.
The degree assignment $\delta(e)$ is obtained from the total $f$-degree of the unbroken connecting
map associated to $e$. 

The $2$-pointed tree $\Gamma_f$ is easily seen to satisfy 
conditions (a)-(e). Condition (c) holds since there are no self connecting maps.
 By construction, $V_3$
is exactly the set of extremal (or edge valence 1) vertices of $\Gamma_f$.
The tree $\Gamma_f$ is invariant as $[f]$ varies in a connected component of the
$\TT$-fixed locus of $\overline{M}_{0,2}(\Hilb_n,d)$. 

Let $G_d$ denote the finite set of $2$-pointed trees of degree $d$.
Let $$\overline{M}_\Gamma \subset \overline{M}_{0,2}(\Hilb_n,d)$$
denote 
the substack of $\TT$-fixed maps corresponding to the tree $\Gamma\in G_d$.
Let
$$\dlang   [\cI_\mu], [\cI_\nu]
  \drang_{d}^{\Gamma}$$
denote the localization contribution of $\overline{M}_\Gamma$.
There is a unique tree in $\Gamma^*\in G$ with
a single edge of degree $d$ corresponding to unbroken maps. Let
$$G^*_d= G_d \setminus \{ \Gamma^* \}.$$
By summing contributions,
\begin{eqnarray*}
\dlang  [\cI_\mu], [\cI_\nu] \drang_{d} & = &
\sum_{\Gamma\in G_d}
\dlang   [\cI_\mu], [\cI_\nu]  \drang_{d}^{\Gamma} \\
& = &
\dlang
 [\cI_\mu], [\cI_\nu] \drang_{d}^{unbroken} + 
\sum_{\Gamma\in G^*_d}
\dlang  [\cI_\mu], [\cI_\nu] \drang_{d}^{\Gamma}
\end{eqnarray*}

\begin{Lemma}\label{ubc}
The $\TT$-equivariant broken contributions vanish.
\end{Lemma}

\begin{proof}
For each $\Gamma \in G^*_d$,  $|E|>1$.
We must show the contribution of each such $\Gamma$ is 0.

Up to automorphisms, the stack $\overline{M}_\Gamma$
factors as a product:
\begin{equation}\label{e23}
\overline{M}_\Gamma = \Big( \prod_{v\in V_1} \overline{M}_v \ \times\  \prod_{e\in E} \overline{M}_e
\Big) \Big/ {\text{Aut}}(\Gamma)
.
\end{equation}
Here, $\overline{M}_v$ denotes the $f$-contracted moduli space of pointed genus 0 
curves associated to $v\in V_1$, and
 $\overline{M}_e$ denotes the moduli space of unbroken $\TT$-fixed maps of degree $\delta(e)$
connecting the $\TT$-fixed points associated to the vertices incident to $e$.
By the virtual localization formula \cite{GrabPan},
\begin{equation}\label{kk12}
\dlang [\cI_\mu], 
[\cI_\nu] \drang_d^\Gamma = \int_{[\overline{M}_\Gamma]^{vir}_r} \frac{\text{ev}_1^*([\cI_\mu]) \cup
\text{ev}_2^*([\cI_\nu])}{e({N}^{vir})},
\end{equation}
where the reduced virtual class on $\overline{M}_\Gamma$ is obtained from the $\TT$-fixed part of the
reduced obstruction
theory and $N^{vir}$ is the virtual normal bundle.

The {\em standard} obstruction theory of $\overline{M}_\Gamma$
is obtained from the $\TT$-fixed part
of the complex 
\begin{equation} \label{xx12}
R^\bullet \pi_*(f^*T_{\Hilb_n})^\vee,
\end{equation}
and the $\TT$-fixed part of the cotangent complex of the Artin stack $L^\bullet({\mathfrak{M}})$.
Here, we follow the notation of Section \ref{rvclass}, see also \cite{GrabPan}.
The 
{\em reduced} obstruction of $\overline{M}_\Gamma$ is obtained by removing a trivial 1-dimensional subobject
from the standard obstruction theory --- see Section \ref{rvclass}.

The normalization
sequence for the universal domain at the breaking curves and breaking nodes relates the 
complex \eqref{xx12} to the corresponding  complexes for each factor in the
product \eqref{e23}.
The normalization sequence on the universal domain  is
$$0 \rightarrow \bigoplus_{v\in V_1} \cO_{C_v} \oplus \bigoplus _{e\in E} \cO_{C_e} \rightarrow
\cO_C \rightarrow \bigoplus_{s\in I} \cO_{s}\rightarrow 0,$$
where $C_v$, $C_e$, are the subcurves associated to $v\in V_1$, $e\in E$, and
$I$ is the set of all incidence points of the subcurves.
After tensoring with the pull-back of $T_{\Hilb_n}$ and taking the derived
$\pi$-push forward to $\overline{M}_\Gamma$, we find the complex \eqref{xx12} differs from the
sum of the corresponding complexes of the factors only by the nodal terms 
$R^\bullet\pi_*( T_{\Hilb_n} \otimes \cO_{s})^\vee$.

The cohomology of the complex associated to a node $s\in I$ is concentrated
in degree 0 and equals the tangent representations 
at $f(s)\in \Hilb^{\TT}_n$. 
The tangent representations at the $\TT$-fixed points $\Hilb_n^{\TT}$
have {\em no} $\TT$-fixed parts. Hence, the $\TT$-fixed part of \eqref{xx12}
is obtained from the sum of $\TT$-fixed parts of corresponding complexes of the
factors \eqref{e23}.

The complex $L^\bullet({\mathfrak M})$
differs from the cotangent complexes of the factors \eqref{e23} by the 
deformation spaces at the nodes $s\in I$ and possible automorphism
factors at the extremal vertices. 
The deformation spaces
at the nodes $s$ have nontrivial $\TT$-weights by definition.
The automorphism factors may differ at the extremal vertices since
the points of $C$ corresponding to $V_3$ may not be marked while
the ends of $\overline{M}_e$ are taken to be marked.
The possible  automorphism factors 
have nontrivial $\TT$-weights (proportional to tangent weights at the
associated $\TT$-fixed points of $\Hilb_n$).
The $\TT$-fixed part of $L^\bullet({\mathfrak{M}})$ is therefore also obtained from the sum of 
cotangent complexes of the factors \eqref{e23}.

We conclude
the standard obstruction theory of
$\overline{M}_\Gamma$ 
is obtained from the sum of the standard
obstruction theories of the factors \eqref{e23}.
For each edge $e$,
the standard obstruction theory of $\overline{M}_e$
admits a trivial 1-dimensional subobject defining the
reduced obstruction theory of $\overline{M}_e$.
The standard obstruction theory of $\overline{M}_\Gamma$
therefore admits a trivial $|E|$-dimensional subobject (compatible, by definition, with
the trivial $1$-dimensional subobject defining by the
reduced obstruction theory of $\overline{M}_\Gamma$). 
Hence, the
reduced obstruction theory of  $\overline{M}_\Gamma$ admits a
trivial $(|E|-1)$-dimensional subobject. If $|E|>1$, 
the reduced virtual class,
$$[\overline{M}_\Gamma]^{vir}_r,$$ 
simply vanishes.
\end{proof}

\subsection{Induction step: II}
\subsubsection{Reduced $3$-point function}
We calculate here the $\TT$-equivariant, reduced, 3-point function
\begin{equation}
\label{ff5567}
  \dlang \left[\cI_{(n)}\right], D, \left[\cI_{(n-1,1)}\right]
\drang_d \, .
\end{equation}
Using the divisor equation and Lemma \ref{ubc},
we write the reduced $3$-point function \eqref{ff5567} as:
\begin{eqnarray*}
\dlang \left[\cI_{(n)}\right], D, \left[\cI_{(n-1,1)}\right]
\drang_d  & =& d \dlang \left[\cI_{(n)}\right], \left[\cI_{(n-1,1)}\right] \drang_d \\
& = &
d\dlang \left[\cI_{(n)}\right], \left[\cI_{(n-1,1)}\right]
\drang_d^{unbroken}.
\end{eqnarray*}


\subsubsection{Unbroken maps}

We must now determine the set of unbroken $\TT$-fixed maps of degree $d$
connecting $\cI_{(n)}$ to $\cI_{(n-1,1)}$.

An {\em unbroken $T$-fixed map} is an  
an unbroken $\TT$-fixed maps which is fixed for the 
{\em full} $T$-action on $\Hilb_n$. If $f$ is
an unbroken $T$-fixed map with an irreducible domain,
then by an analogue of 
 Lemma \ref{444}, we find,
$$w^d_{(n),(n-1,1)} = \frac{-(n-1)t_1+t_2}{d},$$
where $w^d_{(n),(n-1,1)}$ denotes the full $T$-representation.

Since the tangent $T$-weights of $\cI_{(n)}$ lie in a half space, 
the $T$-action on $\Hilb_n$ is isomorphic to a linear $T$-action on a
$T$-invariant affine neighborhood ${\mathbb A}$ of $\cI_{(n)}$.
There is a unique tangent weight of $\cI_{(n)}$ 
proportional to $w^d_{(n),(n-1,1)}$.
The line $L\subset \Hilb_n$,
\begin{equation}
  \label{lineL}
  L: [w_0:w_1] \to (w_0 x^{n-1} + w_1 y, x^n,xy,y^2)\,,
\end{equation}
 is the unique
irreducible, $T$-invariant curve meeting $\cI_{(n)}$
with tangent weight $-(n-1)t_1+t_2$.
Moreover, $L$ 
connects $\cI_{(n)}$ to $\cI_{(n-1,1)}$.

The $d$-fold cover of $L$ is therefore the {\em unique} $T$-fixed map 
of degree $d$ with irreducible domain
connecting $\cI_{(n)}$ and $\cI_{(n-1,1)}$.

\begin{Lemma}
There are no unbroken $T$-fixed maps with reducible domains connecting
$\cI_{(n)}$ to $\cI_{(n-1,1)}$.
\end{Lemma}

\begin{proof}
Let $[f]\in \overline{M}_{0,2}(\Hilb_n,d)$
be a unbroken $T$-fixed map. By definition, $f$ consists of a
chain of non-contracted rational curves
connecting $\cI_{(n)}$ to $\cI_{(n-1,1)}$, 
$${f}: P_1 \cup \ldots \cup P_r\rightarrow \Hilb_n,$$ 
with every node $s_i$ satisfying
$$(t_1+t_2) \ | \ \left(w_{P_i,s_i}+ w_{P_{i+1},s_i}\right),$$
where $w_{P_i,si}$ are the $T$-weights.

We will order partitions by the 
function $\epsilon: {\mathcal P}(n) \rightarrow {\mathbb Z}$,
$$
\epsilon(\lambda)= c(\lambda;1,0)=\sum_{i=1}^{\ell(\lambda)} \binom{\lambda_i}{2}\,,
$$
the $t_1$ coefficient of $c(\lambda;t_1,t_2)$.
By convexity, $\epsilon$ achieves a strict maximum at the partition $(n)$.
The second largest value of $\epsilon$ is achieved uniquely at $(n-1,1)$.

Consider a pair of $T$-fixed maps $h_1$ and $h_2$ 
connecting three points
$$\cI_\mu\ \stackrel{h_1}{\text{---}}  \  \cI_\nu\ \stackrel{h_2}{\text{---}}\  \cI_\xi$$
of $\Hilb_n$.
Assume the maps have irreducible domains $P_1$ and $P_2$ respectively
and the divisibility condition,
$$(t_1+t_2) \ | \ \left( w_{P_1,\nu} + w_{P_2,\nu} \right),$$
is satisfied in the middle.
Let $d_1$ and $d_2$ be the respective degrees of $h_1$ and $h_2$.
By localization,
$$d_1 w_{P_1,\nu} = c(\mu;t_1,t_2)- c(\nu;t_1,t_2),$$
$$d_2 w_{P_2,\nu} = c(\xi;t_1,t_2)- c(\nu;t_1,t_2).$$
The tangent weights
of $\cI_\nu$ are of the form
$$\alpha t_1+ \beta t_2$$
where either 
$\alpha\geq 0, \beta \leq 0$ or $ \alpha\leq 0, \beta \geq 0$.
The $T$-weights $w_{P_1,\nu}$ and $w_{P_2,\nu}$ are proportional to tangent weights
of $\cI_\nu$. Hence, by the divisibility condition and the tangent
weight inequalities,
the $t_1$ coefficients of $w_{P_1,\nu}$ and $w_{P_2,\nu}$
must have opposite signs. We conclude either the condition
$$\epsilon(\cI_\mu) \geq \epsilon(\cI_\nu) \geq  \epsilon(\cI_\xi)$$
or the condition
$$\epsilon(\cI_\mu) \leq \epsilon(\cI_\nu) \leq  \epsilon(\cI_\xi)$$
holds.

For the unbroken $T$-fixed map ${f}$,
the $T$-fixed points $f(s_i)$ must have
$\epsilon$ values lying between
$\epsilon(\cI_{(n)})$ and $\epsilon(\cI_{(n,1)})$.
The latter condition is only possible if, for each node,
$$f(s_i)= \cI_{(n)} \ \ \text{or} \ \ \cI_{(n-1,1)}.$$
Since ${f}$ is reducible, there exists at least
one node. We reach a contradiction since there are
no $T$-fixed maps with irreducible domains connecting
a $T$-fixed point of $\Hilb_n$ to itself.
\end{proof}

The $d$-fold
cover of $L$ is thus the unique unbroken $T$-fixed
map connecting $\cI_{(n)}$ to $\cI_{(n-1,1)}$. 
Since linearized $T$-actions on positive
dimensional varieties must have at least 2 fixed points,
we conclude $dL$ is the unique
unbroken $\TT$-fixed map connecting 
$\cI_{(n)}$ to $\cI_{(n-1,1)}$.

\subsubsection{The contribution of $dL$}
We have proven the equality:
\begin{multline*}
\lang \left[\cI_{(n)}\right],D, 
\left[\cI_{(n-1,1)}\right]\rang_d^{T} \mod (t_1+t_2)^2 = \\
d(t_1+t_2) 
\dlang \left[\cI_{(n)}\right], \left[\cI_{(n-1,1)}\right]
\drang_d^{dL,\TT},
\end{multline*}
where the respective equivariant groups are made explicit
in the notation.
The right side is equal to
$$ d \lang \left[\cI_{(n)}\right], \left[\cI_{(n-1,1)}\right]\rang_{d}^{dL,T}\mod (t_1+t_2)^2 $$
To match the answer of \eqref{JJ}, we will calculate the latter
$T$-equivariant contribution of $dL$.

The contribution is obtained from the $T$-weights of the
representations $$H^0(C,f^*(T_{\Hilb_n})), \ \ H^1(C,f^*(T_{\Hilb_n})),$$
where 
$$f:C \rightarrow L$$ is the unique $T$-fixed unbroken map connecting
$\cI_{(n)}$ and $\cI_{(n-1,1)}$.

For $n>2$,
the restriction of $T_{\Hilb_n}$ to $L$ splits into $T$-equivariant line bundles:
$$
T_{\Hilb_n}\Big|_L=\cO(2)\oplus\cO(-2)\oplus\cO(1)\oplus\cO(-1)\oplus\cO^{2n-4} \,,
$$
where the first summand  is the tangent bundle of $L$. 
The $T$-weights of the trivial part are:
$$
t_1,2t_1,\dots,(n-2)t_1,\quad t_2,t_2-t_1,\dots,t_2-(n-3)t_1 \,.
$$
The corresponding flat deformations of the ideal \eqref{lineL} with 
weight $k t_1$
can be given explicitly by
$$
(w_0 x^{n-1} + w_1 y + \epsilon w_0 x^{n-1-k}, x^n+\epsilon x^{n-k},xy,y^2)\,, 
\quad \epsilon^2=0 \,,
$$
where $k=1,\dots,n-2$. Similarly, 
$$
(w_0 x^{n-1} + w_1 y + \epsilon w_1 x^{k}, x^n,xy+\epsilon x^{k+1},y^2+2\epsilon \delta_{k,0} y)\,, 
\quad \epsilon^2=0 \,,
$$
where $k=0,\dots,n-3$, is a flat deformation with weight $t_2-k t_1$.

The $T$-weights of the nontrivial 
summands are recorded in the following table:

\begin{center}
  \begin{tabular}{c|c|c}
 & $\cI_{(n)}$ & $\cI_{(n-1,1)}$ \\
\hline 
$\cO(2)$ & $t_2-(n-1)\,t_1$ & $(n-1)\,t_1-t_2$ \\
$\cO(-2)$ & $n t_1$ & $2t_2-(n-2)\,t_1$ \\
$\cO(1)$ & $t_2-(n-2)\,t_1$ & $t_1$ \\
$\cO(-1)$ & $(n-1)\,t_1$ & $t_2$ 
  \end{tabular}
\end{center}

For $n=2$, the $T$-equivariant
splitting of $T_{\Hilb_2}$ on $L$ takes a different
form,
$$
T_{\Hilb_2}\Big|_L=\cO(2)\oplus\cO(-2)\oplus\cO\oplus\cO \,.
$$
The $T$-weights are:
\begin{center}
  \begin{tabular}{c|c|c}
 & $\cI_{(2)}$ & $\cI_{(1,1)}$ \\
\hline 
$\cO(2)$ & $t_2-t_1$ & $t_1-t_2$ \\
$\cO(-2)$ & $2 t_1$ & $2t_2$ \\
$\cO$ & $t_2$ & $t_2$ \\
$\cO$ & $t_1$ & $t_1$ 
  \end{tabular}
\end{center}
The weights for the $\cO(1)\oplus \cO(-1)$ summand for $n>2$ are
switched in the $n=2$ case for the $\cO\oplus \cO$ summand.

The $T$-representation $H^0(C,f^*(T_{\Hilb_n}))$ will be shown below to have 
a the single $0$ weight obtained from reparameterization.
As a consequence, $[f]\in \overline{M}_{0,2}(\Hilb_n,d)$
is a nonsingular point of the $T$-fixed locus.
Then, the localization contribution is:
\begin{multline}\label{fred}
d\lang \left[\cI_{(n)}\right], \left[\cI_{(n-1,1)}\right]
\rang_d^{dL} =\\
 \frac{d}{d}\  e\Big(T_{\Hilb_n,\cI_{(n)}}\Big) e\Big(T_{\Hilb_n,\cI_{(n-1,1)}}\Big)
\frac{e\Big(H^1(C,f^*(T_{\Hilb_n}))\Big)}
{e\Big(H^0(C,f^*(T_{\Hilb_n})) - 0\Big)},
\end{multline}
where $e$ denotes the $T$-equivariant Euler class.
The $1/d$ term in front is obtained from the automorphisms of $f$.

We start by calculating the weights of $H^0(C,f^*(T_{\Hilb_n}))$.
The shorthand
$$
\tau = (n-1)t_1-t_2
$$
will be convenient for the formulas.

\begin{enumerate}
\item[$\bullet$]
The weights of $H^0(C,f^*(\cO(2)))$, with the exception of the 0 weight
obtained from reparameterization, 
multiply to 
\begin{multline*}
  \prod_{k=0}^{d-1} \left(t_2-(n-1) t_1 + \frac{k}{d}\, \tau \right)^2 (-1)^{d} 
\equiv\\
 (-1)^{d} 
\left(\frac{t_1 n}{d}\right)^{2d} (d!)^2 \mod (t_1+t_2) \,.
 \end{multline*}
\end{enumerate}
The calculation of $H^0(C,f^*(\cO(1)))$ is separated into two cases:
\begin{enumerate}
\item[$\bullet$]
If $n\notdiv d$, the weights of $H^0(C, f^*(\cO(1)))$ are
\begin{multline}
  \label{H0(1)}
  \prod_{k=0}^{d} \left(t_2-(n-2) t_1 + \frac{k}{d} \, \tau \right) \equiv\\
\left(\frac{t_1 n}{d}\right)^{d+1} 
\frac{\Gamma\left(\frac dn +1\right)}
{\Gamma\left(\frac dn-d\right)}
\mod (t_1+t_2) \,.
\end{multline}

\item[$\bullet$]If $n$ divides $d$, then the factor in \eqref{H0(1)}
corresponding to $k=d-\frac{d}{n}$
equals $\frac1n(t_1+t_2)$ --  reflected by the 
pole of the $\Gamma$-function in the denominator. 
In case $n|d$, the product \eqref{H0(1)}
equals
\begin{equation*}
\frac1n\,  \left(\frac{t_1 n}{d}\right)^{d} 
\frac{\Gamma\left(\frac dn +1\right)}
{\Gamma\left(\frac dn-d+t_1+t_2\right)}
\mod (t_1+t_2)^2 \,.
\end{equation*}
\end{enumerate}
The trivial summands of $T_{\Hilb_n} |_L$ contribute to $H^0(C,f^*(T_{\Hilb_n}))$.
\begin{enumerate}
\item[$\bullet$]
The weights of the trivial summands multiply to
$$
(-1)^{n-2} (n-2)!^2 \,\, t_1^{2n-4} \mod (t_1+t_2)\,,
$$
\end{enumerate}
The computation of the representation $H^0(C, f^*(T_{\Hilb_n}))$ is 
complete, and the 0 weight assertion is verified.

Next, we calculate the weights of $H^1(C,f^*(T_{\Hilb_n}))$. There are
only two summands to consider.

\begin{enumerate}
\item[$\bullet$] The weights of $H^1(C, f^*(\cO(-2)))$ are
\begin{multline*}
  \label{H1(2)}
  \prod_{k=1}^{2d-1} \left(2t_2-(n-2) t_1 + \frac{k}{d} \, \tau \right) \equiv\\
 (-1)^{d-1} (t_1+t_2)\left(\frac{t_1 n}{d}\right)^{2d-2} 
(d-1)!^2
\mod (t_1+t_2)^2 \,.
\end{multline*}
\end{enumerate}
The calculation of $H^1(C,f^*(\cO(-1)))$ is separated into two cases:
\begin{enumerate}
\item[$\bullet$]
If $n\notdiv d$, the weights of $H^1(C, f^*(\cO(-1)))$ are
\begin{equation*}
  \prod_{k=1}^{d-1} \left(t_2+ \frac{k}{d} \, \tau \right) \equiv
(-1)^{d-1} \left(\frac{t_1 n}{d}\right)^{d-1} 
\frac{\Gamma\left(\frac dn\right)}
{\Gamma\left(\frac dn-d+1\right)}
\mod (t_1+t_2) \,.
\end{equation*}
\item[$\bullet$]
When $n|d$, the weights of $H^1(C, f^*(\cO(-1)))$ are
$$
(-1)^{d}\, \frac{n-1}{n} \, \left(\frac{t_1 n}{d}\right)^{d-2} 
\frac{\Gamma\left(\frac dn\right)}
{\Gamma\left(\frac dn-d+1+t_1+t_2\right)}
\mod (t_1+t_2)^2 \,. 
$$
\end{enumerate}

Finally, we require
the Euler classes of $$T_{\Hilb_n, \cI_{(n)}}, \ \ T_{\Hilb_n,\cI_{(n-1,1)}}.$$
The product of the tangent weights at the two points is
$$
\frac{(n!)^4}{(n-1)^2}\, t_1^{4n} \mod (t_1+t_2)\,. 
$$

The contribution of $dL$ is obtained by substituting the weight
calculations in \eqref{fred}. 
We find,  
modulo $(t_1+t_2)^2$,
$$
\lang \left[\cI_{(n)}\right], D, \left[\cI_{(n-1,1)}\right]
\rang_d \equiv 
\begin{cases}
(-1)^{n+d-1}(t_1+t_2) \,\dfrac{t_1^{2n} \, (n!)^2}{n-1} \,, 
& n\notdiv d\,,\\
(-1)^{n+d}(t_1+t_2) \, t_1^{2n} \, (n!)^2\,, 
& n| d\,.
\end{cases}
$$
The generating function for the numbers on the right is 
precisely \eqref{JJ}. The proof of
Theorem \ref{T1} is complete.
\qed

\section{Properties of the quantum ring}

\subsection{Proof of Corollary \ref{cor1}}

The limiting operator,
$$
\lim_{t\to\infty} \frac{1}{t} \, \MM_D\left(q,t,t^{-1}\right) = 
\sum_{k>0} \left( 
\frac{k}{2} \frac{(-q)^k+1}{(-q)^k-1} - \frac{1}{2} \frac{(-q)+1}{(-q)-1}\right)\, \alpha_{-k} \, \alpha_k,
$$
is diagonal with distinct eigenvalues. Hence,
$\MM_D(q,t_1,t_2)$ has distinct eigenvalues for generic
values of the parameters. 

Since the classical ring $H^*_T(\Hilb_n, {\mathbb Q})$ is semisimple after localization,
the quantum ring $QH^*_T(\Hilb_n,{\mathbb Q})$ is also semisimple after localization. 
The idempotents of the quantum ring are
eigenvectors of quantum multiplication by $D$. 
The vector $|1^n\rangle$ represents the unit 
in  $QH^*_T(\Hilb_n, \Q)$. 
Since the unit 
is the sum of 
all idempotents,
the action of $\MM_D$ on $|1^n\rangle$
generates the $n$-eigenvalue subspace of Fock space.
Hence, $D$ generates $QH_T^*(\Hilb_n,\Q)$ after extending
scalars to the field $\Q(q,t_1,t_2)$.
 \qed

\subsection{Multipoint invariants}\label{s_big}

All
$3$-point, genus $0$ $T$-equivariant Gromov-Witten invariants of $\Hilb_n$ in the
Nakajima basis,
$$\lang \lambda, \mu, \nu \rang_{0,3,d}^{\Hilb_n},$$
are determined by Theorem 1 and Corollary 1.
The algorithm below can be used to reconstruct
multipoint genus 0 invariants from $3$-point invariants.

Let $D^{*k}\in QH_T^*(\Hilb_n,\Q)$ denote the $k^{th}$ power of 
$D$ with respect to quantum multiplication. Since the set
$$\{ D^{*k} \}_{0 \leq k \leq |{\mathcal P}(n)|-1},$$ 
spans 
$QH^*_T(\Hilb_n, \Q)$ after extension of scalars,
there is a natural filtration of the quantum ring  
by degree in $D$. 
We will filter the multipoint  
invariants of $\Hilb_n$ by, first, the number of 
insertions $m$ and, second, the minimal degree $k$ in $D$ 
among the insertions.

Since all $3$-point invariants are known, we assume $m\geq 4$.
Since insertions of degree 0 and 
$1$ in $D$ can be removed by the $T$-equivariant fundamental class and divisor 
equations, we assume $k\geq 2$.

Let the following bracket denote a series of $m$-pointed invariants of $\Hilb_n$ of minimal
degree $k$,
$$
\lang D^{*k},\lambda, \mu, \ld \rang^{\Hilb_n} = \sum_{d\geq 0} q^d
\lang D^{*k},\lambda, \mu, \ld \rang^{\Hilb_n}_{0,m,d}.
$$
The dash stands for $m-3$ other insertions. 

Let $\overline{{M}}_{0,4}$ be the moduli space of 4-pointed, genus 0 curves.
Let 
$$(12|34), \ (13|24), \ (14|23) \in \overline{M}_{0,4}$$
denote the three boundary divisors.
Let 
$$
\pi : \overline{{M}}_{0,m+1}(\Hilb_n,d) \to 
\overline{{M}}_{0,4} 
$$
be the $T$-equivariant map obtained by forgetting all the data except for the
first four marking.

Consider the following $(m+1)$-point invariant with domain restriction determined
by $\xi\in \overline{M}_{0,4}$,
\begin{multline*}
\lang D, D^{*(k-1)},\lambda, \mu, \ld  \rang^{\Hilb_n}_\xi = \\
\sum_{d\geq 0} q^d \int_{[\overline{M}_{0,m+1}(\Hilb_n,d)]^{vir}}
\text{ev}_1^*(D) \ \text{ev}_2^*(D^{*(k-1)})\ \text{ev}_3^*(\lambda) \ \text{ev}_4^*(\mu)\ 
( \ld )\ \pi^*([\xi]) .
\end{multline*}
Since points in $\overline{M}_{0,4}$ are cohomologically equivalent, the equality
$$
\lang D, D^{*(k-1)},\lambda, \mu, \ld \rang^{\Hilb_n}_{(12|34)} =
\lang D, D^{*(k-1)},\lambda, \mu, \ld \rang^{\Hilb_n}_{(13|24)}
$$
yields
the WDVV-equation, 
\begin{multline*}
\sum_\nu \lang D,D^{*(k-1)}, \ld,  \nu \rang \, 
\lang \nu^\vee, \lambda, \mu \rang  + 
\lang D, D^{*(k-1)},  \nu \rang \, 
\lang \nu^\vee, \lambda, \mu, \ld \rang = \\
\sum_\nu\lang D,\lambda,  \ld,  \nu \rang \, 
\lang \nu^\vee, D^{*(k-1)},\mu \rang  + 
\lang D, \lambda, \nu \rang \, 
\lang \nu^\vee, D^{*(k-1)},\mu, \ld \rang + \dots. 
\end{multline*}
The summation is over partitions $\nu\in {\mathcal P}(n)$. The
$T$-equivariant Poincare dual of $\nu$ in the Nakajima basis is denoted by $\nu^\vee$. The dots stand for terms 
with nontrivial distribution of the insertions 
(which, therefore, have fewer than $m$ insertions each). 
The superscript $\Hilb_n$ has been dropped from the bracket notation in the WDVV-equation.

By the definition of quantum multiplication, 
$$
\lang D^{*k},\lambda, \mu, \ld \rang^{\Hilb_n} = 
\sum_{\nu}\lang D, D^{*(k-1)}, \nu \rang
\lang \nu^\vee,\lambda,\mu,  \ld  \rang^{\Hilb_n} \,  
$$
All the other 
terms in above WDVV-equation are either 3-point invariants or 
have minimal degree $k-1$ in $D$. 
\qed

\subsection{Relation to the Gromov-Witten theory of $\C^2\times\Pp$}
\label{s_relGW}

We follow here the notation of \cite{jbrp} Section 3.2 for the local 
Gromov-Witten theory of $\C^2 \times \Pp$.

Let $(\Pp,x_{1},\dots
,x_{r})$ be the sphere with $r$ distinct marked points. Let
$$\M (\Pp,\lambda ^{1},\dots ,\lambda ^{r})$$ denote the moduli space of (possibly
disconnected) relative stable maps from genus $h$ curves to $\Pp$ with prescribed ramification
$\lambda ^{i}$ at $x_{i}$.
The prescribed ramification points on the domain are unmarked, 
 and the maps are required to be
nonconstant on all connected components.

The partition function of the local Gromov-Witten theory  
may be defined by:
\begin{multline*}
{\mathsf Z}'_{GW}(\C^2\times \Pp)_{n[\Pp], \lambda^1, \dots, \lambda^r} =\\
\sum _{h\in \Z }
u^{2h-2}
 \int _{[\M (\Pp,\lambda ^{1},\dots ,\lambda ^{r})]^{vir}}e
(-R^{\bullet }\pi _{*}f^{*} (\C^2 \otimes {\mathcal O}_{{\mathbb P}^1})).
\end{multline*}
We will be primarily interested in a shifted generating function,
$$
\mathsf{GW}^*_{n}(\C^2\times \Pp)_{\lambda ^{1},\dots, \lambda ^{r}}= 
(-iu)^{n (2-r)+\sum _{i=1}^{r}\ell (\lambda ^{i})} \ 
{\mathsf Z}'_{GW}(\C^2\times \Pp)_{n[\Pp], \lambda^1, \dots, \lambda^r}.$$

The GW/Hilbert correspondence relates the local theory of $\C^2 \times \Pp$ to
the  multipoint invariants of $\Hilb_n$ with {\em fixed} complex structure
$\xi \in \overline{M}_{0,r}$,
$$
\langle \lambda^1, \dots, \lambda^r \rangle^{\Hilb_n}_{\xi}\,.
$$

\begin{Theorem}\label{gwh} 
After the variable change $e^{iu}=-q$,
$$
\mathsf{GW}^*_{n}(\C^2\times \Pp)_{\lambda^1,\dots, \lambda^r} = (-1)^n \langle \lambda^1, 
\dots, \lambda^r \rangle^{\Hilb_n}_\xi.$$
\end{Theorem}

\begin{proof}
A direct comparison of the formulas of Theorem \ref{T1} of Section \ref{t1} and Theorem 6.5 of \cite{jbrp}
yields the result in case $r=3$ {\em and}  $\lambda^1$ is the 2-cycle $(1^{n-2}2)$. 
A verification shows the degeneration formula
of local Gromov-Witten theory is compatible via the correspondence
with the splitting formula for genus 0 fixed moduli invariants
of $\Hilb_n$.
By Corollary 1, both sides of the correspondence are
canonically determined from the 3-point case with one 2-cycle --- see also the
reconstruction 
result of the Appendix of \cite{jbrp}
\end{proof}

\subsection{The orbifold $(\C^2)^n/S_n$}
Consider the GW/Hilbert correspondence in the 3-point case,
\begin{equation}\label{v234a}
\mathsf{GW}^*_{n}(\C^2\times \Pp)_{\lambda,\mu,\nu} = (-1)^n \langle \lambda, 
\mu,\nu \rangle^{\Hilb_n}.
\end{equation}

The 3-pointed, genus 0, $T$-equivariant Gromov-Witten invariants of
the orbifold  
$(\C^2)^n/S_n$ are easily related to 
$\mathsf{GW}^*_n(\C^2 \times \Pp)_{\lambda,\mu,\nu}$, see
\cite{grby}. 
The Hilbert scheme $\text{Hilb}_n$ is a crepant resolution 
of the (singular) quotient $(\C^2)^n/S_n$.
The equivalence \eqref{v234a}
may be viewed as relating the
$T$-equivariant quantum cohomology
of the quotient {\em orbifold} $(\C^2)^n/S_n$ to the $T$-equivariant quantum cohomology
of the resolution $\text{Hilb}_n$. 

Mathematical conjectures relating the quantum cohomologies of
orbifolds and their crepant resolutions in the
non-equivariant case have been pursued by 
Ruan (motivated by the physical predictions of Vafa and Zaslow). 
Equality  \eqref{v234a} suggests the correspondence also holds
in the equivariant context.

\subsection{Higher genus} \label{higherg}

Localization may be used to compute the 
higher genus Gromov-Witten invariants of $\Hilb_n$. Because 
the $T$-fixed curves are not isolated, the localization 
structure is rather complicated. The
higher genus invariants are expressed as sums over graphs where 
the vertex contributions are Hodge integrals over 
moduli spaces of curves and the edge contributions
are integrals over moduli spaces of 
$T$-fixed curves in $\Hilb_n$. 
The latter can be computed recursively from genus 0 
descendent invariants of $\Hilb_n$. 

We expect the involved localization procedure can be 
conveniently expressed in Givental's formalism 
\cite{giv1,giv2,leep} for higher genus potentials 
for semisimple Frobenius structures. The main 
issue arising in the application of Givental's ideas
is the selection of an $R$-calibration.  We expect the 
standard Bernoulli $R$-calibration used in the $T$-equivariant 
Gromov-Witten theory of toric varieties is appropriate.


\vspace{+10 pt}
\noindent
Department of Mathematics \\
Princeton University \\
Princeton, NJ 08544, USA\\
okounkov@math.princeton.edu \\

\vspace{+10 pt}
\noindent
Department of Mathematics\\
Princeton University\\
Princeton, NJ 08544, USA\\
rahulp@math.princeton.edu

\end{document}